\documentclass{etds}
\usepackage[utf8]{inputenc}
\usepackage{amsmath,amssymb}
\usepackage[british]{babel}
\usepackage{xcolor}
\usepackage{tikz}
\usetikzlibrary{decorations.pathmorphing,positioning,patterns,automata,arrows}
\ifx\UseOption\undefined
\def\UseOption{dessinetikz,dessinempost,bw}
\fi
\usepackage{optional}
\opt{color}{
\colorlet{couleurlien}{green!20!black}
\colorlet{couleurlink}{red!60!black}
}
\opt{bw}{
\colorlet{couleurlien}{black}
\colorlet{couleurlink}{black}
}
\usepackage[pdftex,colorlinks,citecolor=couleurlien,linkcolor=couleurlink]{hyperref}
\newcommand\sshift[1]{\ensuremath{\mathbf{#1}}}
\newcommand\sft{\ensuremath{\sshift{\Sigma}}}
\newcommand\img{\ensuremath{\sshift{\Gamma}}}
\newcommand\src{\ensuremath{\sshift{\Lambda}}}
\newcommand\sofic{\ensuremath{\sshift{X}}}
\newcommand\sofb{\ensuremath{\sshift{Y}}}
\newcommand\sofc{\ensuremath{\sshift{Z}}}

\newcommand\card[1]{\ensuremath{|#1|}}

\newcommand\scover{\ensuremath{\sshift{\Sigma}_R}}
\newcommand\cstar[3]{\ensuremath{d^*_{#3}(#1,#2)}}
\newcommand\cglob[1]{\ensuremath{c^*(#1)}}
\newcommand\preimw[3]{\ensuremath{p_{#3}(#1,#2)}}
\newcommand\limset[1]{\ensuremath{\Omega_{#1}}}

\newcommand\latin[1]{\emph{#1}}

\newcommand\ie{\latin{i.e.}, }
\newcommand\eg{\latin{e.g.}, }
\newcommand\Z{\ensuremath{\mathbb{Z}}}
\newcommand\N{\ensuremath{\mathbb{N}}}

\newcommand\entrop[1]{\ensuremath{h(#1)}}

\newcommand\ooae{$1$-$1$ \emph{a.e.}}

\newcommand\follow[2]{\ensuremath{\mathcal{F}_{#1}^{#2}}}
\newcommand\rer{right-$e$-resolving}
\newcommand\rr{right-resolving}
\newcommand\rerae{\rer{} almost-everywhere}
\newcommand\rrae{\rr{} almost-everywhere}
\newcommand\rcae{right-closing \emph{a.e.}}

\newcommand\rcontae{right-continuing \emph{a.e.}}

\newcommand\pp{\ensuremath{\textrm{Per}}}

\newtheorem{theorem}{Theorem}[section]
\newtheorem{lemma}[theorem]{Lemma}
\newtheorem{prop}[theorem]{Proposition}
\newtheorem{cor}[theorem]{Corollary}
\newtheorem{conj}{Conjecture}

\newtheorem{claim}{Claim}
\newtheorem{question}{Question}

\newcommand\propref[1]{Proposition~\ref{prop:#1}}
\newcommand\claimref[1]{Claim~\ref{claim:#1}}
\newcommand\lemmaref[1]{Lemma~\ref{lemma:#1}}
\newcommand\thmref[1]{Theorem~\ref{thm:#1}}
\newcommand\corref[1]{Corollary~\ref{cor:#1}}
\newcommand\questref[1]{Question~\ref{quest:#1}}

\newcommand\secref[1]{Section~\ref{sec:#1}}
\newcommand\figref[1]{Figure~\ref{fig:#1}}
\newcommand\conjref[1]{Conjecture~\ref{conj:#1}}

\begin{document}

\ETDS{0}{0}{0}{2012}
\runningheads{A. Ballier}{Limit sets of stable CA}

\title{Limit sets of stable Cellular Automata}
\author{Alexis Ballier\affil{1}\footnote{This research has been supported by the FONDECYT Postdoctorado Proyecto
3110088.}
}
\address{\affilnum{1}\ Centro de Modelamiento Matemático\\
        Av. Blanco Encalada 2120 Piso 7\\
        Oficina 710. \\
        Santiago, Chile.\\
\email{aballier@dim.uchile.cl}}

\recd{\today}

\begin{abstract}
        We study limit sets of stable cellular automata standing from a
        symbolic dynamics point of view where they are 
        a special case of sofic shifts admitting a steady epimorphism.
        We prove that there exists a right-closing almost-everywhere steady
        factor map from one irreducible sofic shift onto another one if and only
        if there exists such a map from the domain onto the minimal
        right-resolving cover of the image.
        We define right-continuing almost-everywhere steady maps and prove that
        there exists such a steady map between two sofic shifts if and only if
        there exists a factor map from the domain onto the minimal
        right-resolving cover of the image.

        In terms of cellular automata, this translates into: A sofic shift can
        be the limit set of a stable cellular automaton with a right-closing
        almost-everywhere dynamics onto its limit set if and only if it is the
        factor of a fullshift and there exists a right-closing almost-everywhere
        factor map from the sofic shift onto its minimal right-resolving cover.
        A sofic shift can be the limit set of a stable cellular automaton
        reaching its limit set with a right-continuing almost-everywhere factor
        map if and only if it is the factor of a fullshift and there exists a
        factor map from the sofic shift onto its minimal right-resolving cover.

        Finally, as a consequence of the previous results, we provide a
        characterization of the Almost of Finite Type shifts (AFT) in terms of a
        property of steady maps that have them as range.
\end{abstract}


Cellular automata were introduced by Von Neumann as a model of some biological
processes~\cite{vonneumintroca} and have become a rich model of complex systems:
systems with simple local behavior but complex global evolution. Many different
points of view have been adopted to formalize this complexity, using methods of
combinatorics, topology, ergodic theory, language theory and theory of
computation.
The first, and maybe most known, systematic study of this complex behavior was
performed by S. Wolfram~\cite{Wolfram84c} by doing computer experiments and
then analyzing the observed behavior of the cellular automaton.
From a mathematical point of view, the long-term behavior of a cellular
automaton can be modelled by its dynamics on its limit set: The set of
configurations that can be reached arbitrary late in the evolution of the
automaton.
We can distinguish between two types of cellular automata by their ways of
reaching their limit set, starting from the fullshift, the set of all possible
configurations~\cite{maass95}:
Either the automaton reaches its limit set in finite time, the cellular
automaton is then called \emph{stable}, or it never reaches it and only gets
closer and closer to it, the cellular automaton is then called \emph{unstable}.

In this paper, we are interested in the, \latin{a priori}, simpler case of
stable cellular automata for which it is still an important open problem to
obtain a characterization of their limit sets~\cite[Section~$16$]{boyleop}.
Stable cellular automata can be modelled in terms of symbolic
dynamics~\cite{marcuslind,maass95}: They are a special case of the steady factor
maps of Barth and Dykstra~\cite{barth2007}.
Basic remarks yield necessary conditions for a subshift to be the limit set of a
stable cellular automaton: This is what A. Maass called property
(H)~\cite{maass95}. A. Maass then proved that these necessary conditions are
also sufficient for a large class of sofic shifts: The almost of finite type
(AFT \cite{marcusdefaft}) shifts~\cite[Theorem~$4.8$]{maass95}.
Albeit not exactly stated as such in A. Maass' paper, his methods for
constructing limit sets of stable cellular automata are to obtain a weak
conjugacy between two subshifts (constructing factor maps from each subshift
onto the other one) and then if one can prove that one subshift is a
stable limit set of cellular automata then the other one is automatically
also a stable limit set~\cite[Lemma~$4.1$]{stableunstable}.
As a consequence of Boyle's extension
lemma~\cite[Lemma~$2.4$]{boylelowentropfact}, a subshift of finite type (SFT)
having property (H) is the stable limit set of a cellular
automaton~\cite[Theorem~$3.2$]{maass95}.
Moreover, all the methods we know for constructing stable cellular automata go
through a weak conjugacy with an SFT, this is what led us to the following
conjecture that we restate here:
\begin{conj}{\cite[Conjecture~$1$]{stableunstable}}
        \label{conj:slswccov}
The limit set of any stable cellular automaton is weakly conjugate to an SFT.
\end{conj}

After fixing the definitions in \secref{defs}, where we adopt the
vocabulary from symbolic dynamics~\cite{marcuslind}, we prove the basic results
that we will use along the rest of the paper.
In \secref{rcae}, we prove that a sofic shift is the stable limit
set of a cellular automaton with a right-closing almost-everywhere dynamics on
its limit set if and only if the sofic shift has property (H) and there exists a
right-closing almost-everywhere factor map from the sofic shift onto its minimal
right-resolving cover.
In \secref{rcontaeca}, we prove that a cellular automata attains its limit set
by a right-continuing almost-everywhere factor map if and only if its limit set
factors onto its minimal right-resolving cover.
By similar methods, we provide in \secref{newcaraft} a characterization of the
almost of finite type (AFT) shifts of B.~Marcus~\cite{marcusdefaft} in terms of
the range of a special class of steady maps and characterize AFT stable limit
sets of cellular automata as those that can be attained by a left and
right-continuing almost-everywhere cellular automaton.

Each of these three sections (\ref{sec:rcae}, \ref{sec:rcontaeca} and
\ref{sec:newcaraft}) are organized in the same way and each of them provides a
characterization in terms of steady maps (Theorems~\ref{thm:caracsoficrcaestd},
\ref{thm:factmapextrcaessiwccov} and \ref{thm:factmapextrcaessiwccovaft}
respectively). One direction of each of these characterizations always makes use
of an extension theorem for sliding block codes: these are, respectively,
Boyle's extension lemma~\cite[Lemma~$2.4$]{boylelowentropfact}, its refinement
by Boyle and Tuncel~\cite[Theorem~$5.3$]{inftoonecodes} and yet another
refinement by Jung~\cite[Theorem~$4.5$]{bicontcodes}.
Hence, Sections~\ref{sec:rcae}, \ref{sec:rcontaeca} and
\ref{sec:newcaraft} are organized in a somewhat chronological order of the
results they are based on.

\section{Definitions and basic results}

\label{sec:defs}

Let $A$ be a finite set, called the \emph{alphabet} embedded with the discrete
topology.
Consider $A^{\Z}$ as the \emph{fullshift} over $A$ embedded with the product
topology.
For $i\in\Z$ and $x\in{}A^{\Z}$, denote by $x_i$ the value of $x$ at position
$i$.
A metric for the topology of $A^{\Z}$ can be defined for example as
$d(x,y)=2^{-\min\{|i|, x_i\neq{}y_i\}}$.

\subsection*{Words and languages}
A \emph{word} over $A$ is an element of $A^*=\cup_{n\in\N}A^n$.
Denote by $|w|$ the \emph{length} of the word $w$, \ie such that
$w\in{}A^{|w|}$.
For $i<j\in\Z$ and $x\in{}A^{\Z}$, denote by $x_{[i;j]}$ the word
$x_ix_{i+1}\ldots{}x_j\in{}A^*$.
We say that a word $w$ \emph{appears} in $x\in{}A^{\Z}$ at position $i$ if
$x_{[i;i+|w|-1]}=w$.
For a subset $\sofic$ of $A^{\Z}$, we can define the \emph{language} of $\sofic$
as the set of words that appear in some element of $\sofic$:
$\mathcal{L}(\sofic)=\left\{w\in{}A^*, \exists x\in\sofic, \exists i\in\Z,
x_{[i;i+|w|-1]}=w\right\}$. 
To ease notations we denote, for $x\in{}A^{\Z}$, we denote $\mathcal{L}(x)$ for
$\mathcal{L}(\left\{x\right\})$.
When $w\in\mathcal{L}(\sofic)$, we say that $w$ is an \emph{$\sofic$-word}.
Denote by $\mathcal{L}_n(\sofic)$ the set of words of length $n$ appearing in
$\sofic$, \ie $\mathcal{L}_n(\sofic)=\mathcal{L}(\sofic)\cap{}A^n$.

\subsection*{Shift and subshifts}
Define the \emph{shift} $\sigma:A^{\Z}\to{}A^{\Z}$ as $\sigma(x)_i=x_{i+1}$.
$\sigma$ is bijective, thus induces a $\Z$-action on the fullshift $A^{\Z}$.
A subset $\sofic$ of $A^{\Z}$ is said to be \emph{shift-invariant} if
$\sigma(\sofic)=\sofic$.
A \emph{subshift} of $A^{\Z}$ is a \emph{closed and
shift-invariant} subset of $A^{\Z}$.

\subsection*{Transitive and asymptotic configurations}
For a subshift $\sofic$, a configuration $x\in\sofic$ is said to be
\emph{left-transitive} in $\sofic$ if $\mathcal{O}_-(x)=\left\{\sigma^i(x),
i\leq{}0\right\}$ is dense in $\sofic$.
It is \emph{right-transitive} in $\sofic$ if
$\mathcal{O}_+(x)=\left\{\sigma^i(x), i\geq{}0\right\}$ is dense in $\sofic$.
Two configurations $x,y\in{}A^{\Z}$ are said to be \emph{left-asymptotic} if
there exists $n\in\Z$ such that for all $i\leq{}n$, $x_i=y_i$. They are
\emph{right-asymptotic} if there exists $n\in\Z$ such that for all $i\geq{}n$,
$x_i=y_i$.

\subsection*{Forbidden words}
It is well known that a subshift can also be defined by a set of \emph{forbidden
words} $\mathcal{F}\subseteq{}A^*$: $\sofic$ is a subshift of $A^{\Z}$ if and
only if there exists $\mathcal{F}\subseteq{}A^*$ such that
$\sofic=\left\{x\in{}A^{\Z}, \forall w\in\mathcal{F},
w\not\in\mathcal{L}(x)\right\}$.
The above $\mathcal{F}$ can be always chosen as
$A^*\setminus{}\mathcal{L}(\sofic)$.
When such an $\mathcal{F}$ can be chosen \emph{finite} we say that $\sofic$ is a
\emph{subshift of finite type}, SFT in short.
If the length of the longest word of such a finite $\mathcal{F}$ is not greater
than $2$ then it is said to be a \emph{one-step SFT}.

\subsection*{Factor maps}
Let $\src$ and $\img$ be subshifts. A map $f:\src\to{}\img$ is
\emph{shift-commuting} if $\sigma\circ{}f=f\circ\sigma$.
A continuous, shift-commuting and onto map $f:\src\to{}\img$ is called a
\emph{factor map}.
A bijective factor map is called a \emph{conjugacy}.
If there exist factor maps $\pi:\sofic\to\sofb$ and $\varphi:\sofb\to\sofic$
then we say that the subshifts $\sofic$ and $\sofb$ are \emph{weakly conjugate}.
A \emph{sofic shift} is the image of an SFT by a factor map.
It is clear that a subshift conjugate to an SFT or a sofic shift is itself,
respectively, an SFT or a sofic shift.
If $\pi:\sft\to\sofic$ is a factor map from an SFT onto a sofic shift then
$(\sft,\pi)$ is called a \emph{cover} of $\sofic$.

\subsection*{Sliding block codes}
For $D$ a finite subset of $\Z$, a \emph{block map} on $D$ is a function
$g:A^{D}\to{}B$ where $A$ and $B$ are finite sets.
$g$ defines a \emph{sliding block code} $f:A^{\Z}\to{}B^{\Z}$ by
$f(x)_i=g(x_{|D+i})$. 
When $D=\left\{0\right\}$, $f$ is said to be \emph{one-block}.
By the Curtis-Hedlund-Lyndon theorem~\cite{hedlund69}, sliding block codes
between $A^{\Z}$ and $B^{\Z}$ are exactly the continuous and shift-commuting
maps between those spaces.
Among other things, this implies that a bijective sliding block code (\ie a
conjugacy) has a sliding block code inverse.

\subsection*{Magic words}
Let $f:\src\to\img$ be a one-block factor map.
Let $m=m_1\ldots{}m_k$ be a $\img$-word.
For $0<i\leq{}k$, let $\cstar{m}{i}{f}$ be the number of different symbols we can
see at position $i$ in a $f-$pre-image of $m$, that is:
$\cstar{m}{i}{f}=\card{\preimw{m}{i}{f}}$ where $\preimw{m}{i}{f}=\left\{c_i,
f(c)=m\right\}$.
Denote by $\cglob{f}$ the minimum of $\cstar{m}{i}{f}$ over all $i\in\N$ and
all $\img$-words $m$.
A word $m$ such that $\cstar{m}{i}{f}=\cglob{f}$ is a called a \emph{magic word}
for $f$ at coordinate $i$.

The following property of magic words will help in understanding better the
notions we use in this paper:
\begin{prop}[Mainly {\cite[Corollary~$9.1.10$]{marcuslind}}]
        \label{prop:mwallpreim}
        Let $f:\src\to\img$ be a one-block factor map and $m$ a magic word for
        $f$ at coordinate $i$. For any $\img$-word of the form $vmw$ (that is,
        an extension of $m$) and any symbol $c\in\preimw{m}{i}{f}$ there exists
        an $f-$pre-image $VMW$ of $vmw$ such that $M_i=c$.
\end{prop}

\proc{Proof.}
        Note that any $\src$-word $VMW$ that is an $f$-pre-image of $vmw$ is
        such that $M_i\in\preimw{m}{i}{f}$, that is
        $\preimw{vmw}{i}{f}\subseteq{}\preimw{m}{i}{f}$.
        But since $w$ is magic at coordinate $i$, we have
        $\card{\preimw{m}{i}{f}} = \cstar{m}{i}{f} \leq{} \cstar{vmw}{i}{f} =
        \card{\preimw{vmw}{i}{f}}$, hence $\preimw{vmw}{i}{f}=\preimw{m}{i}{f}$.
\ep
\medbreak

\subsection*{Entropy}
For a subshift $\sofic$, one can define its \emph{entropy}, which roughly
speaking represents the exponential growth rate of its language:
\[
\entrop{\sofic}=\lim_{n\to\infty}\frac{\log\card{\mathcal{L}_n(\sofic)}}{n}
\]
For example, if $\sofic$ and $\sofb$ are conjugate subshifts then
they have the same entropy. If $\sofic$ factors onto $\sofb$ then
$\entrop{\sofb}\leq{}\entrop{\sofic}$, that is, entropy does not increase
\latin{via} factor maps.

\subsection*{Irreducibility and mixing}
A subshift $\sofic$ is said to be \emph{irreducible} if for any two
configurations $x,y\in\sofic$, there exists $N\in\N$ and $z\in\sofic$ such that
$z_i=x_i$ for $i\leq{}0$ and $z_i=y_i$ for $i\geq{}N$.
It is well known that if $\sofic$ is sofic then there exists such an $N$ that
does not depend on the configurations $x$ and $y$.
$\sofic$ is said to be \emph{mixing} if there exists $N\in\N$ such that for any 
$k\geq{}N$ and any two configurations $x,y\in\sofic$ there exists 
$z\in\sofic$ such that $z_i=x_i$ for $i\leq{}0$ and $z_i=y_i$ for $i\geq{}k$.
A factor of an irreducible subshift is itself irreducible
and a factor of a mixing subshift is also mixing.

\subsection*{Fiber product}
\begin{figure}
        \begin{center}
            \opt{dessinetikz}{\begin{tikzpicture}[scale=.6,->,>=stealth']
\tikzstyle{every state}=[draw=none];

\node[state] (f)                           {$\sshift{F}$};
\node[state] (fp)   [below      =of f  ]   {$\sshift{F}'$};
\node[state] (t1)   [below right=of fp ]   {};
\node[state] (t2)   [below left =of fp ]   {};
\node[state] (sft1) [      left =of t2 ]   {$\sft_1$};
\node[state] (sft2) [      right=of t1 ]   {$\sft_2$};
\node[state] (sof)  [below right=of t2 ]  {$\sofic$};

\path[->] (f)   edge                 node[above left ]     {$\rho_1$}       (sft1)
          (f)   edge                 node[above right]     {$\rho_2$}       (sft2)
          (fp)  edge                 node[below right]     {$\rho_{1|\sshift{F}'}$}       (sft1)
          (fp)  edge                 node[below left ]     {$\rho_{2|\sshift{F}'}$}       (sft2)
          (sft1) edge                node[below left ]     {$\pi_1$}       (sof)
          (sft2) edge                node[below right]     {$\pi_2$}       (sof)
                ;
\path[white]     (fp) edge node[black] {\rotatebox[origin=c]{90}{$\subseteq$}} (f);
\end{tikzpicture}}\opt{pdftikz}{\includegraphics{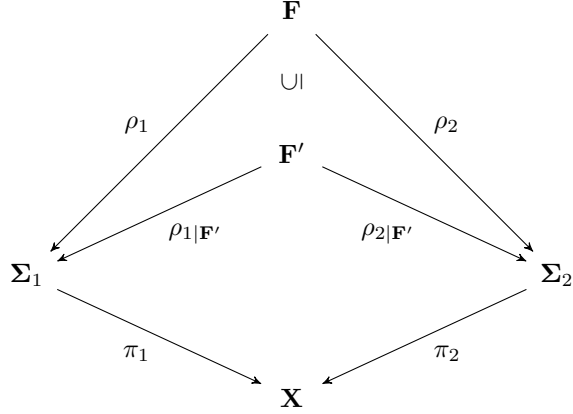}}
        \end{center}
        \caption{A commutative diagram of the fiber product and full fiber product of $(\sft_1,\pi_1)$ and $(\sft_2,\pi_2)$.}
        \label{fig:fiberpbase}
\end{figure}

A classical construction from symbolic
dynamics~\cite[Definition~$8.3.2$]{marcuslind} is the \emph{fiber product} of
two covers of the same sofic shift: Let $(\sft_1,\pi_1)$ and $(\sft_2,\pi_2)$ be
covers of the same sofic shift $\sofic$.
We define the \emph{full fiber product} $\sshift{F}$ of  $(\sft_1,\pi_1)$ and
$(\sft_2,\pi_2)$ as: $\sshift{F}=\left\{(x_1,x_2), x_1\in\sft_1, x_2\in\sft_2,
\pi_1(x_1)=\pi_2(x_2)\right\}$. $\sshift{F}$ comes with canonical projections:
$\rho_1:\sshift{F}\to\sft_1, \rho_1(x_1,x_2)=x_1$ and
$\rho_2:\sshift{F}\to\sft_2, \rho_1(x_1,x_2)=x_2$.
Usually, $\rho_1$ inherits the properties of $\pi_2$ and $\rho_2$ those of
$\pi_1$~\cite[Proposition~$8.3.3$]{marcuslind}; we will state precisely what
this means when we will need it.
Since both $\sft_1$ and $\sft_2$ are SFTs, so is $\sshift{F}$.
If $\sft_1$ and $\sft_2$ are irreducible, then $\sshift{F}$ is not necessarily
irreducible, however, it contains a unique irreducible component of maximal
entropy $\sshift{F}'$ and the restrictions of $\rho_1$ and $\rho_2$ to
$\sshift{F}'$ remain surjective. $\sshift{F}'$ is called the \emph{fiber
product} of $(\sft_1,\pi_1)$ and $(\sft_2,\pi_2)$.
The situation after all those definitions is depicted on \figref{fiberpbase}.

\subsection*{Right-closing and resolving}
A factor map $f:\src\to\img$ between subshifts is said to be
\emph{right-closing} if it never collapses two left-asymptotic points.
That is:
        \[
                \forall x,y\in\src, \forall i\leq{}0, x_i=y_i,
                f(x)=f(y)\Rightarrow x=y
        \]
$f$ is \emph{right-closing almost-everywhere} if we only require the
above to hold for left-transitive $x$ and $y$.
        If $f$ is a one-block map, it is said to be \emph{right-resolving} if
        whenever $ab$ and $ac$ are two-letters words in $\src$ then $f(b)=f(c)$
        implies $b=c$.

\subsection*{Minimal right-resolving cover}
Among the covers of an irreducible sofic shift $\sofic$, there is one of
particular interest: the \emph{minimal right-resolving cover}, or \emph{Fischer
cover}~\cite{fischersofic} $(\scover,\pi_R)$. $\pi_R:\scover\to\sofic$ is
a one-block right-resolving factor map and if $f:\sft\to\sofic$ is a
right-closing factor map then there exists $\varphi:\sft\to\scover$ such that
$f=\pi_R\circ\varphi$~\cite[Proposition~$4$]{caracaftmincover}.

\subsection*{Periodic points}
A configuration $x$ is said to be \emph{periodic} if there exists an integer
$i>1$ such that $\sigma^i(x)=x$. The \emph{period} of $x$ is the smallest such
$i$. We denote by $\pp(\sofic{})$ the \emph{set of periodic points} of the
subshift $\sofic$.
If a subshift $\sofic$ factors onto a subshift $\sofb$ then for every periodic
point $x$ of $\sofic$ there exists a periodic point $y$ of $\sofb$ whose period
divides the period of $x$ (take $y$ to be the image of $x$ by the factor map).
We denote this relation $\pp(\sofic)\to\pp(\sofb)$.
It turns out that this trivial necessary condition on periodic points is also
sufficient for the existence of a factor map between two irreducible SFTs of
unequal entropy~\cite{boylelowentropfact}.
A periodic point $x$ is represented as a finite word $w$, whose length is the
period of $x$, repeated infinitely: $x={}^\infty{}w^\infty$.
Following \cite{boylelowentropfact}, for an irreducible sofic shift $\sofic$
with minimal right-resolving cover $(\scover,\pi_R)$, we say that such an $x$ is
a \emph{receptive periodic point} if there exist magic words for $\pi_R$: $m_1$
and $m_2$ such that for every $n\geq{}1$, $m_1w^nm_2$ is an $\sofic$-word.
If $\sofic$ is SFT then any periodic point is receptive because $\pi_R$ is a
conjugacy. Following \cite{maass95} we say that a configuration $x$ is a
\emph{receptive fixed point} if it is a receptive periodic point of period $1$.
As remarked at the end of section $2$ in \cite{maass95}, a factor map between
irreducible sofic shifts maps receptive fixed points to receptive fixed points.

\subsection*{Right-resolving almost everywhere}
        If $\src$ is an irreducible sofic shift with minimal right-resolving
        cover $(\scover,\pi_R)$ and $f$ is one-block, we say that $f$ is
        \emph{right-resolving almost-everywhere} if $f$ is right-closing
        almost-everywhere and $f\circ\pi_R:\scover\to\img$ is right-resolving.

\subsection*{Right-continuing and \rer{}}
        A factor map $f:\src\to\img$ between sofic shifts is said to be
        \emph{right-continuing}~\cite{inftoonecodes}
        if for any $x$ in $\src$ and $y$ in $\img$ such
        that $f(x)$ and $y$ are left-asymptotic, there exists $x'$
        left-asymptotic to $x$ in $\src$ such that $f(x)=y$.
        If there exists an integer $n$ such that for any $x\in\src$ and
        $y\in\img$ such that $f(x)_{(-\infty;n]}=y_{(-\infty;n]}$ then there
        exists $x'\in\src$ such that $x'_{(-\infty;0]}=x_{(-\infty;0]}$ and
        $f(x')=y$ then $f$ is said to be \emph{right-continuing with retract
        $n$}.
        If $f$ is right-continuing with retract $0$ then it is said to be
        \emph{\rer{}}.
        As before, we define \emph{right-continuing almost-everywhere} and
        \emph{\rer{} almost-everywhere} when we only require the above to hold
        for a left-transitive $y$ but impose the existence of a retract.
        Indeed, \propref{mwallpreim} implies that \emph{any} factor map with an
        SFT domain is right (and left) continuing almost-everywhere, but without
        retract.
        We will usually not append ``with a retract'' when talking about
        right-continuing almost-everywhere factor maps and consider the
        existence of the retract to be part of the definition of
        right-continuing almost-everywhere.
        By \cite[Proposition~$5.1$, (iii)$\Rightarrow$(i)]{inftoonecodes}, if
        $f$ is right-continuing (resp. \rcontae{}) with a retract then there
        exists a conjugacy $\Theta$ and $f'$ such that $f=f'\circ\Theta$ and
        $f'$ is \rer{} (resp. \rerae{}).
        Also remark that in the definition of \rcae{} we imposed $x$ to be
        left-transitive while in the definition of \rcontae{} we impose $y$ to
        be left-transitive: it is simply a matter of historical definitions,
        \rcae{} has, to our knowledge, always been defined as such while we
        could impose $y$ to be left-transitive in the definition of \rcae{}
        since for a finite-to-one factor map $f$, $f(x)$ is left-transitive if
        and only if $x$ is by a slight modification of
        \cite[Lemma~$9.1.13$]{marcuslind}.

\subsection*{Links between right-closing and right-continuing with a retract}

Right-continuing shall be seen as the dual of right-closing and \rer{} the dual
of right-resolving.
One may remark that the above definition of \rer{} is more intricate than the
original one for SFTs in \cite{inftoonecodes} and than its right-resolving
dual; they are equivalent when $\src$ is SFT but differ when it is merely sofic:
With the original definition we may have \rer{} factor
maps over sofic shifts which are not right-continuing~\cite{contcodesyu}.
The above definition avoids this problem and is equivalent to the original one
for SFTs by \cite[Proposition~$5.1$]{inftoonecodes}.

While the right-continuing image of an SFT is an SFT~\cite{contcodesyu} (or
~\cite[Proposition~$2.1$]{boyleputnam} for the finite-to-one case, or even
\cite[Proposition~$8.2.2$]{marcuslind}), a
right-closing factor map from an SFT is right-continuing
almost-everywhere~\cite[Lemma~$2.5$]{boyleputnam}.
A right-closing almost-everywhere factor map with SFT domain is right-closing
(everywhere)~\cite[Proposition~$4.10$]{bmt}.
Therefore for a finite-to-one $f$, we may ask whether right-closing
almost-everywhere is equivalent to right-continuing almost-everywhere.

\begin{prop}[Mainly {\cite[Lemma~$2.5$]{boyleputnam}}]
        \label{prop:rcaercontae}
        A factor map $f:\src\to\img$ between irreducible sofic
        shifts which is right-closing almost-everywhere is also
        right-continuing almost-everywhere (with a retract).
        If $f$ is right-resolving almost-everywhere then it is \rerae{}.
\end{prop}

\proc{Proof.}
        Let $f:\src\to\img$ be right-closing almost-everywhere.
        Let $(\scover,\pi_R)$ be the minimal right-resolving cover of $\src$.
        $f\circ\pi_R:\scover\to\img$ is right-closing
        almost-everywhere~\cite[Proposition~$4.11$]{bmt} and thus
        right-closing~\cite[Proposition~$4.10$]{bmt}.
        By \cite[Lemma~$2.5$]{boyleputnam}, $f\circ\pi_R$ is right-continuing
        almost-everywhere with a retract.
        Let $x$ be left-transitive in $\src$ and $\tilde{x}$ its (left-transitive)
        pre-image in $\scover$.
        Let $y$ be left-asymptotic to $f(x)$ in $\img$. Since $f\circ\pi_R$ is
        right-continuing almost-everywhere, find $\tilde{x}'$ in $\scover$,
        left-asymptotic to $\tilde{x}$ such that $f\circ\pi_R(\tilde{x}')=y$.
        $x'=\pi_R(\tilde{x}')$ is the $x'$ we were looking for.
        The right-resolving case follows similarly to
        \cite[Lemma~$2.5$]{boyleputnam}: $x'_{i+1}$ and $x_{i+1}$ are uniquely
        determined by $y_{i+1}$ and, respectively, $x'_i$ and $x_i$; since
        $y_i=f(x_i)$ for $i\leq{}0$ and $x'_i=x_i$ for $i\leq{}-n$, then
        $x'_i=x_i$ for $i\leq{}0$.
\ep
\medbreak

The converse of \propref{rcaercontae} holds when $f$ is finite-to-one:

\begin{prop}        
        \label{prop:rraef21rcae}
        If a finite-to-one factor map $f:\src\to\img$ between irreducible sofic
        shifts is \rcontae{} (with a retract) then it is right-closing
        almost-everywhere.
\end{prop}

\proc{Proof.}
        Up to a conjugacy we can assume that $f$ is \rerae{}.
        Let $f:\src\to\img$ be \rerae{} and suppose it is not right-closing
        almost-everywhere.
        Since, from \propref{rcaercontae}, $\pi_R$ is \rerae{},
        $f\circ\pi_R:\scover\to\img$ is also \rerae{}.
        By \cite[Proposition~$4.11$]{bmt}, $f\circ\pi_R$ is right-closing
        almost-everywhere if and only if $f$ is.
        Therefore, by considering $f\circ\pi_R$ we can assume that $\src$ is a
        one-step SFT. Since $f$ is finite-to-one, by
        \cite[Proposition~$9.1.7$]{marcuslind}, we may assume that $f$ has a
        magic symbol $b$.

        Let $x$ and $y$ be two left-transitive left-asymptotic configurations of
        $\src$ such that $f(x)=f(y)=z$.
        Without loss of generality, suppose $x_i=y_i$ for all $i<0$ and
        $x_0\neq{}y_0$.
        By irreducibility of $\img$, let $z'$ be a right-transitive configuration
        of $\img$ such that for all $i\leq{}0$, $z_i=z'_i$.
        Since $f$ is \rerae{}, let $x'$ and $y'$ be configurations of $\src$
        such that $f(x')=f(y')=z'$ and for all $i\leq{}0$, $x_i=x'_i$ and
        $y_i=y'_i$.
        Since $z'$ is bi-transitive, let $j<0$ and $k>0$ be such that
        $z'_j=z'_k=b$.
        $x'_{[j;k]}$ and $y'_{[j;k]}$ are two pre-images of a word starting and
        ending by the magic symbol $b$, therefore they are mutually separated by
        \cite[Proposition~$9.1.9$]{marcuslind}: they are either equal or differ
        in every coordinate.
        However, $x'_{-1}=x_{-1}=y_{-1}=y'_{-1}$ and $x'_0=x_0\neq{}y_0=y'_0$
        and since $-1$ and $0$ are in the interval $[j;k]$,  $x'_{[j;k]}$ and
        $y'_{[j;k]}$ cannot be mutually separated, a contradiction.
\ep
\medbreak

Note that we cannot remove the hypothesis on the retract in
\propref{rraef21rcae}: Otherwise since \propref{mwallpreim} implies that 
\emph{any} factor map is right-continuing almost-everywhere without retract, any
finite-to-one factor map from an SFT would be right-closing almost-everywhere
and thus right-closing by \cite[Proposition~$4.10$]{bmt}, however there exist
finite-to-one factor maps between SFTs that are not right-closing.

The following proposition shall
be seen as the dual of \cite[Proposition~$4.12$]{bmt} which states that a
\rcae{} factor map from an SFT onto a sofic shift is right-closing everywhere:

\begin{prop}
    \label{prop:reraeimgstfrer}
    A \rerae{} factor map $f:\src\to\img$, where $\src$ is an
    irreducible sofic shift and $\img$ is an irreducible SFT, is \rer{}
    (everywhere).
\end{prop}

\proc{Proof.}
    Let $x\in\src$ and $y\in\img$ be such that $f(x)_i=y_i$ for all $i\leq{}0$.
    For an integer $n$, find by irreducibility of $\src$ a left-transitive
    configuration $x^n\in\src$ such that $x^n_i=x_i$ for all $i\geq{}-n$.
    Define $y^n$ such that $y^n_i=y_i$ for $i\geq{}0$ and
    $y^n_i=f(x^n)_i$ for $i<0$. For $n$ sufficiently big, $y^n$ belongs to
    $\img$ since it is SFT.
    $y^n$ is left-transitive because $x^n$ is and $f$ is onto.

    Since $f$ is \rerae{}, find $z^n$ such that $f(z^n)=y^n$ and $z^n_i=x^n_i$
    for $i\leq{}0$.
    By compactness of $\src$, we may assume w.l.o.g. that $z^n$ converges to
    $z\in\src$. $y^n$ clearly converges to $y$, thus, by continuity of $f$,
    $f(z)=y$.
    Moreover, for all $i\leq{}0$, $z_i=x_i$, thus $f$ is \rer{} everywhere.
\ep
\medbreak

Again, in \propref{reraeimgstfrer}, \rerae{} can be replaced by right-continuing
almost-everywhere with a retract and we get a right-continuing with a retract
factor in the conclusion.
Without the retract hypothesis, it may be possible that the $z^n$ we find agrees
with $x^n$ only at positions $i<-n$ so that its limit may not be left-asymptotic
to $x$ at all.

\begin{prop}
        \label{prop:rcontaedecompfaible}
        Let $\Phi:\sofic\to\sofb$ and $\Psi:\sofb\to\sofc$ be factor maps
        between irreducible sofic shifts. If $\Psi\circ\Phi:\sofic\to\sofc$ is
        right-continuing almost-everywhere with a retract then so is $\Psi$.
\end{prop}

\proc{Proof.}
        Let $N$ be the retract of $\Psi\circ\Phi$.
        Assume without loss of generality that $\Psi$ and $\Phi$ are both
        one-block.
        Let $y\in\sofb$ be left-transitive and $z\in\sofc$ be such that
        $\Psi(y)_i=z_i$ for $i\leq{}N$.
        Let $x\in\sofic$ be a left-transitive pre-image of $y$.
        By our hypothesis, there exists $x'\in\sofic$ such that $x'_i=x_i$ for
        $i\leq{}0$ and $\Psi\circ\Phi(x')=y$. Let $y'=\Phi(x')$.
        Since $\Phi$ is one-block, $y'_i=y_i$ for $i\leq{}0$.
        And $\Psi(y')=z$, thus $\Psi$ is right-continuing almost-everywhere with
        retract $N$.
\ep
\medbreak

\subsection*{Followers}

        Let $\sft$ be a one-step SFT and $f:\sft\to\img$ be a one-block factor
        map onto a sofic shift.
        For a letter $a$ of the alphabet of $\sft$, denote $\follow{f}{\sft}(a)$
        the \emph{$f-$follower} of $a$:
        $\follow{f}{\sft}(a)=\left\{f(aw),aw\in\mathcal{L}(\sft)\right\}$.
        $(\sft,f)$ is said to be \emph{follower
        separated}~\cite[Definition~$3.3.7$]{marcuslind} if for any letters $a$
        and $b$ of the alphabet of $\sft$, if
        $\follow{f}{\sft}(a)=\follow{f}{\sft}(b)$ then $a=b$.

        Remark that the minimal right-resolving cover of a sofic shift is always
        follower separated~\cite[Proposition~$3.3.9$]{marcuslind}, it is
        actually the only (up to conjugacy) cover that is both follower
        separated and right-resolving. It is also well known that any factor map
        from an SFT can be decomposed through a follower separated factor map
        with SFT domain:

\begin{lemma}[{\cite[Section~$3.3$]{marcuslind}} or the remarks before
        {\cite[Proposition~$1.2$]{trowlift}}]
        \label{lemma:folsep}
        Let $\sft$ be a one-step SFT and $f:\sft\to\img$ a one-block factor map
        onto a sofic shift. There exists a one-step SFT $\tilde{\sft}$ and
        one-block factor maps $\varphi:\sft\to\tilde{\sft}$ and
        $\pi:\tilde{\sft}\to\img$ such that $f=\pi\circ\varphi$ and
        $(\tilde{\sft},\pi)$ is follower-separated.
\end{lemma}

In~\cite[Proposition~$1.2$]{trowlift}, it is proved, in addition, that when $f$
is finite-to-one then $\varphi$ can be chosen right-resolving.

\subsection*{Cellular automata and their limit set}
A \emph{cellular automaton} is a sliding block code $f:A^{\Z}\to{}A^{\Z}$, \ie
an endomorphism of a fullshift.
The \emph{limit set} of a cellular automaton $f$ is denoted by $\limset{f}$ and
is defined by: \[ \limset{f}=\bigcap_{n\in\N} f^n(A^{\Z}) \]
One can prove by a simple compactness argument that $\limset{f}$ is precisely
the set of configurations that have $f^n$-pre-images for any integer $n$.
A cellular automaton is said to be \emph{stable}~\cite{maass95} if there exists
an integer $N$ such that $\limset{f}=f^N(A^{\Z})$; in other words the cellular
automaton reaches its limit set in finitely many steps.
$\limset{f}$ is always closed and shift-invariant, hence a subshift.
Since $A^{\Z}$ is a mixing sofic shift with a receptive fixed point, so is
$\limset{f}=f^N(A^{\Z})$ when $f$ is a stable cellular automaton. These are the
necessary conditions for being a stable limit set of cellular automaton that A.
Maass called property (H) in \cite{maass95}.
In \cite{maass95} he also proved that these conditions happen
to be sufficient for almost of finite type shifts (see \secref{newcaraft}
for the definition and more details on these sofic shifts).
However, it is an important open problem to get a characterization of
such subshifts that can occur as limit sets of cellular automata in the general
case, even for the, \latin{a priori}, simpler case where we assume the cellular
automata to be stable~\cite[Section~$16$]{boyleop}.

\subsection*{Steady maps}

\begin{figure}
        \begin{center}
            \opt{dessinetikz}{\begin{tikzpicture}[scale=.6,->,>=stealth']
\tikzstyle{every state}=[draw=none];

\node[state] (sft)                    {$\sft$};
\node[state] (xs)   [below =of sft]   {$\src$};
\node[state] (xt)   [right =of xs ]   {$\img$};

\path[->] (sft) edge                 node[above]     {$f$}       (xt)
          (xs)  edge                 node[below]     {$f$}       (xt)
                ;
\end{tikzpicture}}\opt{pdftikz}{\includegraphics{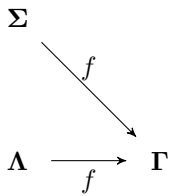}}
        \end{center}
        \caption{A steady factor map.}
        \label{fig:steady}
\end{figure}

Let $\src$ and $\img$ be irreducible sofic shifts.
A factor map $f:\src\to\img$ is said to be \emph{steady}~\cite{barth2007}, or
$\sft$-steady, if there exists an SFT $\sft$ containing $\src$ such that $f$ is
well defined on $\sft$ and $f(\sft)=f(\src)=\img$.
The diagram of a steady map is represented on \figref{steady}.

Steady maps provide a good formalism for stable limit set of cellular automata:
if $f:A^{\Z}\to{}A^{\Z}$ is a stable cellular automaton, then there exists $N$
such that $f^N(A^{\Z})=\limset{f}$. By definition of $\limset{f}$, we have
$f^N(\limset{f})=\limset{f}$, thus, $f^N:\limset{f}\to\limset{f}$ is a steady
epimorphism of $\limset{f}$.
This means that stable cellular automata are a special kind of steady maps
between irreducible sofic shifts. In the rest of the paper we will focus on
steady maps and state the results we obtain for them as theorems while their
implications on stable limit set of cellular automata will be corollaries, even
if characterizing limit sets of stable cellular automata is what motivated our
study.

If $f:\src\to\img$ is a $\sft$-steady factor map, then we say that $f$ is a
\emph{right-closing almost-everywhere steady map} if $f:\src\to\img$ is, in
addition, right-closing almost-everywhere. $f$ is a \emph{right-continuing
almost-everywhere steady map} if $f:\sft\to\img$ is right-continuing
almost-everywhere (with a retract). Note that the domain on which we consider
$f$ differs between the two definitions. The former class of steady maps is
studied in \secref{rcae} and the latter in \secref{rcontaeca}.

For a stable cellular automaton $f:A^{\Z}\to{}A^{\Z}$, let $N$ be an integer
such that $f^N(A^{\Z})=\limset{f}$. We say that $f$ is a 
\emph{right-closing almost-everywhere cellular automaton} when $f^N$ is a
right-closing almost-everywhere steady map, $f$ is a \emph{right-continuing
almost-everywhere stable cellular automaton} when $f^N$ is a right-continuing
almost-everywhere steady map.

\section{Right-closing almost everywhere steady maps}

\label{sec:rcae}

In this section we study right-closing almost-everywhere steady maps. We prove
that such maps can always be decomposed through the minimal right-resolving
cover of its range (\lemmaref{rcaethendecomp}) and then characterize sofic
shifts between which there can exist such a factor (\thmref{caracsoficrcaestd})
so that we get a characterization of sofic shifts that are stable limit set of
cellular automata with a \rcae{} dynamics on its limit set
(\corref{rcaeca}).

\begin{lemma}
        \label{lemma:rcaethendecomp}
    If $f:\src\to\img$ is a right-closing almost-everywhere $\sft$-steady factor
    map then $f$ can be decomposed through the minimal right-resolving cover of
    $\img$ by a right-continuing factor map with a retract, \ie there
    exists $\varphi:\src\to\scover$ such that $f=\pi_R\circ\varphi$ and
    $\varphi$ is right-continuing with a retract.
\end{lemma}

\proc{Proof.}
        Without loss of generality (up to a conjugacy), we can assume that $f$
        is one-block and that $\sft$ is the one-step SFT approximation of
        $\src$ such that $f(\sft)=f(\src)=\img$.

        Let $\varphi$ and $\pi$ be the factor maps given by applying
        Lemma~\ref{lemma:folsep} to $(\sft,f)$.
        Let $\tilde{\src}=\varphi(\src)$ and $\tilde{\sft}=\varphi(\sft)$.
        This is summarized in the following diagram:

\begin{center}
        \opt{dessinetikz}{\begin{tikzpicture}[scale=.6,->,>=stealth']
\tikzstyle{every state}=[draw=none];

\node[state] (sft)                         {$\sft$};
\node[state] (xtd)  [      right=of sft]   {$\tilde{\sft}$};
\node[state] (xt)   [below      =of xtd]   {$\tilde{\src}$};
\node[state] (xs)   [      left =of xt ]   {$\src$};
\node[state] (img)  [      right=of xt ]   {$\img$};

\path[->] (sft) edge                 node[above]       {$\varphi$}       (xtd)
          (xs)  edge                 node[below]       {$\varphi$}       (xt)
          (xtd) edge                 node[above right] {$\pi$}           (img)
          (xt)  edge                 node[below]       {$\pi$}           (img)
                ;

\path[white]     (xs) edge node[black] {\rotatebox[origin=c]{90}{$\subseteq$}} (sft);
\path[white]     (xtd) edge node[black] {\rotatebox[origin=c]{90}{$\subseteq$}} (xt);
\end{tikzpicture}
}\opt{pdftikz}{\includegraphics{dessins/diagdecomp.pdf}}
\end{center}

By \cite[Proposition~$4.11$]{bmt}, both $\varphi:\src\to\tilde{\src}$ and
$\pi:\tilde{\src}\to\img$ are right-closing almost-everywhere.
Without loss of generality, we can assume that they are right-resolving
almost-everywhere.
By \propref{rcaercontae}, $\pi:\tilde{\src}\to\img$ is \rerae{}.

\begin{claim}
        \label{claim:pirr}
        $\pi:\tilde{\src}\to\img$ is right-resolving (everywhere).
\end{claim}

\proc{Proof of \claimref{pirr}.}
        Suppose that we have $ab_1$ and $ab_2$ allowed in $\tilde{\src}$ such
        that $\pi(b_1)=\pi(b_2)=b$.
        Since $(\tilde{\sft},\pi)$ is follower separated,
        $\follow{\pi}{\tilde{\sft}}(b_1)\neq\follow{\pi}{\tilde{\sft}}(b_2)$.
        Without loss of generality, let $w$ be a word such that $b_1w$ is
        allowed in $\tilde\sft$ and
        $\pi(w)\in\follow{\pi}{\tilde{\sft}}(b_1)\setminus\follow{\pi}{\tilde{\sft}}(b_2)$.
        Let $z$ be a right-transitive word of $\tilde\sft$ starting by $w$ and 
        $xab_2y$ be a left-transitive configuration of $\tilde{\src}$;
        $xab_1z$ is a valid configuration of $\tilde{\sft}$ since it is a
        one-step SFT, thus,
        $\pi(xab_1z)=\pi(xa)\pi(z)$ is a configuration of $\img$.
        Since $\pi:\tilde{\src}\to\img$ is \rerae{},
        there exists $z'$ such that $xab_2z'$ is a configuration of
        $\tilde{\src}$ and
        $\pi(xab_2z')=\pi(xa)\pi(z)=\pi(xa)\pi(w)\pi(z)_{[|w|;\infty)}$.
        However, we assumed that
        $\pi(w)\not\in\follow{\pi}{\tilde{\sft}}(b_2)$,
        a contradiction.
        We therefore conclude that such $ab_1$ and $ab_2$ cannot exist and thus
        that $\pi$ is right-resolving everywhere.
\ep
\medbreak

        Since $\pi:\tilde{\src}\to\img$ is right-closing and
        $\pi(\tilde{\sft})=\img$, by \cite[proof of
        Proposition~$4.2$]{maass95}, $\tilde{\src}$ is SFT.
        By \cite[Proposition~$4$]{caracaftmincover}, $\pi$ being a
        right-closing cover of $\img$ can be decomposed through $\scover$, which
        proves the existence of the decomposition: $f=\pi_R\circ\varphi$.

        Now, since $f$ is right-closing almost-everywhere, so is $\varphi$ by
        \cite[Proposition~$4.11$]{bmt}. By taking a conjugacy, we may assume
        that $\varphi$ is \rrae{}, but $\varphi$ maps $\src$ onto $\scover$
        which is SFT, thus by \propref{reraeimgstfrer}, $\varphi$ is \rer{}
        everywhere. By unwinding the conjugacy we took at the beginning, we get
        that the original $\varphi$ is right-continuing with a retract.
\ep
\medbreak

\begin{theorem}
        \label{thm:caracsoficrcaestd}
        Let $\src$ and $\img$ be irreducible sofic shifts of equal entropy.
        There exists a \rcae{} steady factor map from $\src$ onto $\img$ if and
        only if there exists a right-continuing (with a retract) factor map from
        $\src$ onto the minimal right-resolving cover of $\img$.
\end{theorem}

\proc{Proof.}
    $\Rightarrow$: \lemmaref{rcaethendecomp}.

    $\Leftarrow$: Let $\varphi:\src\to\scover$ be a right-continuing factor map
    from $\src$ onto the minimal right-resolving cover of $\img$.
    $\src$ and $\scover$ having the same entropy, $\varphi$ is finite-to-one and
    thus \rcae{} by \propref{rraef21rcae}.
    $\pi_R\circ\varphi$ is thus a right-closing almost-everywhere factor map
    from $\src$ onto $\img$ by \cite[Proposition~$4.11$]{bmt}.
    
    It remains to prove that $\pi_R\circ\varphi$ is steady.
    Let $\sft$ be an irreducible SFT containing $\src$ such that
    $\pp(\sft)\to\pp(\scover)$. Note that such an SFT $\sft$ always exists by
    \eg \cite[Lemma~$4.1$]{ashley}.
    Now we can apply Boyle's extension
    lemma~\cite[Lemma~$2.4$]{boylelowentropfact} to extend $\varphi$ to
    $\tilde{\varphi}:\sft\to\scover$. $f=\pi_R\circ\tilde{\varphi}$ is therefore
    the desired \rcae{} steady factor map.
\ep
\medbreak

\begin{cor}
        \label{cor:rcaeca}
        A subshift is the stable limit set of a cellular automaton which
        has a \rcae{} dynamics on its limit set if and only if it is a factor of
        a fullshift\footnote{Being a factor of a fullshift is equivalent to the property (H) of A.
        Maass in \cite{maass95} by Boyle's lower entropy factor theorem for sofic
shifts~\cite[Theorem~$3.3$]{boylelowentropfact}.} and factors by a
right-continuing factor map onto its minimal right-resolving cover.
\end{cor}

\section{Right-continuing almost-everywhere steady maps}
\label{sec:rcontaeca}

In this section we study right-continuing almost-everywhere steady maps.
It is organized the same way as \secref{rcae}: We prove that such maps can
always be decomposed through the minimal right-resolving cover of its range
(\lemmaref{rcontaedoncdecomp}) and then characterize sofic shifts between which
there can exist such a factor (\thmref{factmapextrcaessiwccov})
so that we get a characterization of sofic shifts that are stable limit set of
cellular automata that attain their limit set with a right-continuing
almost-everywhere factor map (\corref{slrcaewc}).

\begin{lemma}
        \label{lemma:rcontaedoncdecomp}
        A factor map $f:\sft\to\img$ that is right-continuing almost-everywhere
        with retract $N$ from an irreducible SFT $\sft$ onto a sofic shift $\img$
        can be decomposed through the minimal right-resolving cover
        $(\scover{},\pi_R)$ of $\img$.
\end{lemma}

\proc{Proof.}
        Without loss of generality we can assume that $\sft$ is a one-step SFT
        and $f$ is one-block.

        Let $\sshift{F}$ be the irreducible component of maximal entropy of the
        fiber-product of $(\sft,f)$ and $(\scover,\pi_R)$. Let
        $\rho_1:\sshift{F}\to\sft$ and $\rho_2:\sshift{F}\to\scover$ be the
        canonical projections.
        By \cite[Proposition~$5.1$]{trowdecomp} $f$ can be decomposed through
        $(\scover,\pi_R)$ if and only if $\rho_1$ is a conjugacy.

        Suppose $\rho_1$ is not a conjugacy: Let
        $(x,y^1)$ and $(x,y^2)$ be two configurations of $\sshift{F}$ that have
        the same $\rho_1$-image.
        Without loss of generality, suppose $y^1_N\neq{}y^2_N$.
        Since $y^1_N$ and $y^2_N$ are symbols of $\scover$ and $\pi_R$ is
        follower-separated, we can assume that there exists $w$ such that
        $y^1_Nw$ is allowed in $\scover$ and no $y^2_Nw'$, where
        $\pi_R(w')=\pi_R(w)$, is allowed.

        By irreducibility of $\sshift{F}$, find a left-transitive configuration
        $(x',y')\in\sshift{F}$ such that $(x'_i,y'_i)=(x_i,y^1_i)$ for
        $i\geq{}0$.
        Let $y''\in\scover$ be such that $y''_i=y'_i$ for $i\leq{}N$ and
        $y''_{[N+1;N+|w|]}=w$.
        Since $f$ is right-continuing almost-everywhere with retract $N$, there
        exists $x''$ such that $f(x'')=\pi_R(y'')$ and $x''_i=x'_i$ for
        $i\leq{}0$. Hence, $x''_0=x_0$.
        
        Since $x''_0=x_0$ and $\sft$ is one-step, there exists $z\in\sft$ such
        that $z_i=x_i$ for $i\leq{}0$ and $z_i=x''_i$ for $i>0$.
        Since $\pi_R$ is right-resolving, by classical fiber-product
        arguments~\cite[Proposition~$8.3.3$]{marcuslind}, $\rho_1$ is also
        right-resolving.
        By \cite[Proposition~$2.1(5)$]{boyleputnam}, since both $\sshift{F}$ and
        $\sft$ are irreducible, $\rho_1$ is \rer{}, and
        thus there exists $y^3$ such that $(z,y^3)\in\sshift{F}$ and
        $y^3_i=y^2_i$ for $i\leq{}0$.
        Now, $\pi_R(y^3)=f(z)$ and $f(z)_i=\pi_R(y^2)_i$ for $i\leq{}N$, hence
        since $\rho_1$ is right-resolving, $y^3_i=y^2_i$ for $i\leq{}N$.
        Let $w'=y^3_{[N+1;N+|w|]}$.
        Since $f(z)=\pi_R(y^3)$ and $f(z)_{[N+1;N+|w|]}=\pi_R(w)$ we have $w'$
        such that $y^2_Nw'$ is allowed in $\scover$ and $\pi_R(w')=\pi_R(w)$, a
        contradiction.
\ep
\medbreak

\begin{theorem}
        \label{thm:factmapextrcaessiwccov}
        Let $\src$ and $\img$ be irreducible sofic shifts.
        There exists a right-continuing almost-everywhere steady factor map
        $f:\src\to\img$ if and only if $\src$ factors onto $\scover$, the
        minimal right-resolving cover of $\img$.
\end{theorem}

\proc{Proof.}
        $\Rightarrow$:
        If $f:\sft\to\img$ is right-continuing almost-everywhere with a retract
        then by  \lemmaref{rcontaedoncdecomp}, there exists
        $\varphi:\sft\to\scover$, where $(\scover,\pi_R)$ is the minimal
        right-resolving cover of $\img$ such that $f=\pi_R\circ\varphi$.
        Then, $\varphi(\src)\subseteq\scover$. But then
        $f(\src)=\pi_R(\varphi(\src))=\img$, and $\pi_R$ is finite-to-one so
        that $\varphi(\src)$ and $\img$ have the same entropy. Hence,
        $\varphi(\src)$ and $\scover$ also have the same entropy, and since
        $\varphi(\src)$ and $\scover$ are both irreducible sofic shifts, they
        are actually equal: $\varphi:\src\to\scover$ is onto.

        $\Leftarrow:$ Let $\varphi:\img\to\scover$ be the factor map from $\img$
        onto $\scover$ and $\pi_R:\scover\to\img$ be the minimal right-resolving
        cover of $\img$.
        Let $\sft$ be an irreducible SFT containing $\src$ such that
        $\pp(\sft)\to\pp(\scover)$. Note that such an SFT $\sft$ always exists
        by \cite[Lemma~$4.1$]{ashley}.

        By \cite[Theorem~$5.3$]{inftoonecodes}, $\varphi$ can be extended to a
        right-continuing factor map $\tilde{\varphi}:\sft\to\scover$.
        $\tilde{\varphi}$ has a retract and $\pi_R:\scover\to\img$ is \rerae{}
        by \propref{rcaercontae} or \cite[Lemma~$2.5$]{boyleputnam}.
        Let $N$ be the retract of $\tilde{\varphi}:\sft\to\scover$.

        Let $f=\pi_R\circ\tilde{\varphi}$. It is clear that $f:\src\to\img$ is
        a $\sft$-steady factor map.
        We claim that $f:\sft\to\img$ is right-continuing almost-everywhere with
        retract $N$: Let $x\in\sft$ and $y\in\img$ be a left-transitive
        configuration such that $f(x)_i=y_i$ for $i\leq{}N$.
        Let $x'=\tilde{\varphi}(x)\in\scover$ and $y'\in\scover$ a
        $\pi_R$-pre-image of $y$: We have $\pi_R(x')=f(x)$ and $\pi_R(y')=y$.
        Since $f(x)$ and $y$ are left-transitive and $\pi_R$ is \ooae{}
        right-resolving, it is clear that $x'_i=y'_i$ for $i\leq{}N$.
        Now, since $\tilde{\varphi}$ is right-continuing with retract $N$, there
        exists $z\in\sft$ such that $\tilde{\varphi}(z)=y'$ and $z_i=x_i$ for
        $i\leq{}0$.
        Then $f(z)=\pi_R(\tilde{\varphi}(z))=\pi_R(y')=y$ and $f$ is indeed
        right-continuing almost-everywhere with retract $N$.
\ep
\medbreak

\begin{cor}
        \label{cor:slrcaewc}
        A subshift is the stable limit set of a right-continuing
        almost-everywhere cellular automaton if and only if it is
        a factor of a fullshift and is
        weakly conjugate to its minimal right-resolving cover.
\end{cor}

\section{AFT shifts}

\label{sec:newcaraft}

In this section, we continue with the same methods we used in the previous two
sections to obtain a characterization of AFT shifts by means of the type of
steady maps that have them as range (\thmref{factmapextrcaessiwccovaft}).
As is usually the case, the situation is much simpler in the AFT case and the
conclusions can be strengthened.

        An irreducible sofic shift $\img$ is said to be \emph{AFT}, for
        \emph{Almost of Finite Type}, if its minimal right-resolving cover
        $(\scover,\pi_R)$ is also left-closing~\cite{marcusdefaft}.

\begin{lemma}
        \label{lemma:bicontaeaft}
        Let $f:\sft\to\img$ be a factor map from an irreducible SFT $\sft$ onto
        a sofic shift $\img$.
        If $f$ is right and left-continuing almost-everywhere with a bi-retract
        then $\img$ is AFT.
\end{lemma}

\proc{Proof.}
        By \lemmaref{rcontaedoncdecomp}, $f$ can be decomposed through
        $(\scover,\pi_R)$.
        Then, by \propref{rcontaedecompfaible} (and its analogous replacing
        right-continuing by left-continuing), $\pi_R$ is right and
        left-continuing almost-everywhere with a bi-retract.
        Then, by \propref{rraef21rcae}, $\pi_R$ being finite-to-one is right
        and left-closing almost everywhere.
        Finally, by \cite[Proposition~$4.10$]{bmt}, since $\scover$ is SFT,
        $\pi_R$ is both right and left-closing everywhere and thus $\img$ is
        AFT.
\ep
\medbreak

\begin{theorem}
        \label{thm:factmapextrcaessiwccovaft}
        Let $\src$ and $\img$ be irreducible sofic shifts.
        There exists a left and right-continuing almost-everywhere steady factor
        map $f:\src\to\img$ if and only if $\src$ factors onto the minimal
        right-resolving cover of $\img$ and $\img$ is AFT.
\end{theorem}

\proc{Proof.}
        $\Rightarrow$: \lemmaref{bicontaeaft}.

        $\Leftarrow$: 
        Let $\varphi:\src\to\scover$ be the factor map in the hypothesis from
        $\src$ onto the minimal right-resolving cover of $\img$.
        Let $\sft$ be an irreducible SFT containing $\src$ such that
        $\pp(\sft)\to\pp(\scover)$. As before, such an SFT $\sft$ always exists
        by \cite[Lemma~$4.1$]{ashley}.
        By \cite[Theorem~$4.5$]{bicontcodes}, $\varphi$ can be extended to
        a bi-continuing factor map $\tilde{\varphi}:\sft\to\scover$. Let
        $f'=\pi_R\circ\tilde{\varphi}:\sft\to\img$.
        Since $\img$ is AFT, $\pi_R$ is right and left-continuing
        almost-everywhere with a bi-retract, therefore so is $f'$ as the
        composition of two such maps.
\ep
\medbreak

\begin{cor}
\label{cor:bicontaft}
        A subshift is the stable limit set of a left and right-continuing
        almost-everywhere cellular automaton if and only if it is
        a factor of a fullshift and is AFT.
\end{cor}

\proc{Proof.}
        The $\Rightarrow$ direction is clear from
        \thmref{factmapextrcaessiwccovaft}. The $\Leftarrow$ direction requires
        a bit more work:
        By \cite[Lemma~$4.9$]{maass95}, the minimal right-resolving cover of an
        AFT shift that is a factor of a fullshift has a fixed point for
        $\sigma$.
        By \cite[Corollary~$1.3$]{ashley}, there exists a factor map from the
        AFT onto its minimal right-resolving cover since the condition on
        periodic points is fulfilled by the existence of a fixed point in the
        cover. Then we can apply
        \thmref{factmapextrcaessiwccovaft}.
\ep
\medbreak

\section{Conclusions and questions}

\label{sec:questconcl}

We characterized the existence of certain steady maps between irreducible sofic
shifts by the existence of certain factors onto the minimal right-resolving
covers of the image, thus providing characterization of the limit sets of
certain stable cellular automata.
The most annoying problem is that we do not know if there exist limit sets of
stable cellular automata that cannot be reached by (possibly another) cellular
automaton with such properties, meaning \conjref{slswccov} remains a conjecture.

\begin{figure}
        \begin{center}
            \opt{dessinetikz}{\begin{tikzpicture}[scale=.6,shorten >=1pt,node distance=3cm,on grid,auto]

\node[state] (a)    {};
\node[state] (b)  [right =of a]  {};
\node[state] (c)  [above =of a]  {};

\path[->] 
          (a)  edge [bend right]    node[right]    {$1$}    (c)
               edge [bend left]    node[above] {$2$} (b)
          (c)  edge [bend right]    node[left ]    {$1$}    (a)
               edge [bend left]    node[above right] {$2,3$} (b)
          (b)  edge [bend left]    node[below ]    {$4$}    (a)
                ;
\end{tikzpicture}}\opt{pdftikz}{\includegraphics{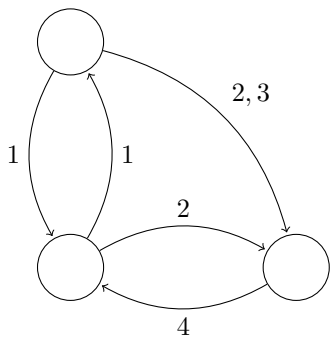}}
        \end{center}
        \caption{A (non-AFT) mixing sofic shift with a receptive fixed point but
        whose cover has no fixed point.}
\label{fig:covnonfix}
\end{figure}

We may note that there exist sofic shifts that have receptive fixed points (and
thus are factor of a fullshift by Boyle's lower entropy factors theorem for
sofic shifts~\cite[Theorem~$3.3$]{boylelowentropfact})
but whose minimal right-resolving cover does not have a fixed point as depicted
on \figref{covnonfix}:
On \figref{covnonfix}, ${}^{\infty}1^{\infty}$ is a fixed point.
It is also a receptive fixed point: $41^*2$ is always a valid word and $2$ and
$4$ are magic.
Therefore, by \corref{slrcaewc}, this subshift cannot be obtained as the stable
limit set of a cellular automaton that is right-continuing almost-everywhere.
We do not know if this subshift is a stable limit set of cellular automata:

\begin{question}
        Is the subshift depicted on \figref{covnonfix} a stable limit set of
        cellular automaton ?
\end{question}

Note that the minimal left-resolving cover of this sofic shift has a fixed
point, this is to keep the example simple; it is left to the reader to modify it
so that neither the minimal left nor right resolving covers have a fixed point.
By \cite[Lemma~$4.9$]{maass95}, an AFT shift which is a factor of a fullshift
always has a minimal right-resolving cover with a fixed point. This example
shows that it is not the case in general.

We may also note that there exist sofic shifts that have \emph{no} equal entropy
SFT factor~\cite{coredimconst}.
The example provided in \cite{coredimconst}(Example~$2.3$) is even worse: It has
a receptive fixed point and is such that its minimal right and left-resolving
covers both have fixed points, showing that the periodic points obstruction is
not the only one. We do not know if \cite[Example~$2.3$]{coredimconst} can be
the stable limit set of a cellular automaton.

Remark that all the stable limit sets of cellular automata constructed in
\cite{maass95} and \cite{stableunstable} have a \rcae{} dynamics on their limit
sets.
Therefore, by \corref{rcaeca} and \corref{slrcaewc}, these subshifts can be
obtained by a right-continuing almost-everywhere and right-closing
almost-everywhere cellular automaton.
There, obviously, exist stable cellular automata that do not have a \rcae{}
dynamics on their limit set: consider any non right-closing surjective CA.
However, it may be possible that its limit set can also be attained by a
cellular automaton with such a property:

\begin{question}
        \label{quest:limrcae}
        If $\sshift{X}$ is the stable limit set of a cellular automaton, is it
        the stable limit set of a cellular automaton which is right-closing
        almost-everywhere on its limit set?
        Of a \rcontae{} cellular automaton ?
\end{question}

Since a subshift that is weakly conjugate to the stable limit set of a CA is
itself the stable limit set of (another) CA~\cite[Lemma~$4.1$]{stableunstable}, we
may weaken \questref{limrcae} to the following:

\begin{question}
        \label{quest:wcrcae}
        Is a stable limit set of CA weakly conjugate to a subshift for which
        such a right-closing or right-continuing almost-everywhere cellular
        automaton exists ?
\end{question}

Remark that \questref{wcrcae} is equivalent to \conjref{slswccov}: 
If every stable limit set of CA is weakly conjugate to an SFT then since
basically any onto endomorphism of an SFT (\eg the identity which is
\rcae{}) can be obtained as the dynamics of a stable CA on an SFT limit
set~\cite{maass95} and can also be attained by a right-continuing factor map.
Conversely, by \thmref{caracsoficrcaestd} or \thmref{factmapextrcaessiwccov}, if
there exists such a CA then its limit set is weakly conjugate to an SFT.

A way to construct stable limit sets of CA (and actually, the only method we
know) is to prove they are weakly conjugate
to an SFT. Moreover, in the constructions, this SFT is always the minimal
right-resolving cover of the limit set.
We may ask if this is the only way to do it with this technique:

\begin{question}
        \label{quest:wccovrcae}
        Let $\src$ be an irreducible sofic shift and $(\scover,\pi_R)$ its
        minimal right-resolving cover. If $\src$ factors onto $\scover$, does it
        factor onto $\scover$ with a right-closing almost-everywhere factor map?
\end{question}

Remark also that by a theorem of J. Ashley~\cite{ashleyres}, two SFTs are weakly
conjugate if and only if they are weakly conjugate by right-closing factor maps.
We cannot require the same for sofic shifts since a right-closing factor map
with SFT range has a SFT domain, but we may ask if it remains true by replacing
right-closing by \rcae{}:

\begin{question}
    \label{quest:wcwcrcae}
    Are two weakly conjugate irreducible sofic shifts also weakly conjugate by
    right-closing almost-everywhere factor maps?
\end{question}
   
\questref{wccovrcae} is a special case of \questref{wcwcrcae} since if a sofic
shift factors onto its minimal right-resolving cover then they are weakly
conjugate.

If the answer to \questref{wcwcrcae} is positive, meaning J. Ashley
results~\cite{ashleyres} can be extended to sofic shifts, then as a consequence
of Boyle's extension lemma~\cite[Lemma~$2.4$]{boylelowentropfact}, we can
construct a \rcae{} steady epimorphism of any sofic shift that is weakly
conjugate to an SFT and thus by \thmref{caracsoficrcaestd} this sofic shift is
weakly conjugate to its minimal right-resolving cover. This would mean that we
can replace SFT by ``the minimal right-resolving cover of the sofic shift'' in
\conjref{slswccov}.
Remark that the answer to \questref{wcwcrcae} is positive for AFT shifts as soon
as the trivial condition on periodic points is satisfied, for which the proof is
short enough to include it here:
\begin{prop}
        \label{prop:wcrcaewc}
        If $\sofic$ and $\sofb$ are two weakly conjugate AFT shifts, with their
        respective minimal right-resolving covers
        $(\Sigma_{\sofic},\pi_{\sofic})$ and $(\Sigma_{\sofb},\pi_{\sofb})$ such
        that $\pp(\sofic)\to\pp(\Sigma_{\sofb})$ and
        $\pp(\sofb)\to\pp(\Sigma_{\sofic})$, then there exist right-closing
        almost-everywhere factor maps from $\sofic$ onto $\sofb$ and from
        $\sofb$ onto $\sofic$.
\end{prop}

\proc{Proof.}
Let $\varphi:\sofic\to\sofb$ be a factor map from $\sofic$ onto $\sofb$.
By \cite[Theorem~$9$]{caracaftmincover}, $\varphi\circ\pi_{\sofic}$ can be
decomposed through $\pi_{\sofb}$ so that $\Sigma_{\sofic}$ factors
$\Sigma_{\sofb}$.
By the same reasoning, $\Sigma_{\sofb}$ factors onto $\Sigma_{\sofic}$ so that
they are weakly isomorphic SFTs; by \cite[Corollary~$1.2$]{ashleyres}, their
dimension groups are isomorphic.

Since $\pp(\sofic)\to\pp(\Sigma_{\sofb})$, by \cite[Theorem~$1.2$]{ashley},
$\sofic$ factors onto $\Sigma_{\sofb}$ by a \rcae{} factor map, and thus also
factors onto $\sofb$ by a \rcae{} factor map.
By the same reasoning, $\sofb$ factors onto $\sofic$ by a \rcae{} factor map.
\ep
\medbreak

Remark that by \cite[Lemma~$4.9$]{maass95} the condition on periodic points hold
if the AFT shifts have a receptive fixed point and therefore the answer to
\questref{wcwcrcae} is positive for AFT shifts with a receptive fixed point.

\bibliographystyle{plain}
\bibliography{biblio}

\end{document}